\documentclass[a4paper, english, 12pt]{amsart}

\usepackage[english]{babel}
\usepackage[latin1]{inputenc}
\usepackage{latexsym}
\usepackage{amsmath,amsthm,amsfonts,amssymb}
\usepackage{dsfont,mathrsfs}

\newtheorem{theo}{Theorem}
\newtheorem{prop}[theo]{Proposition}
\newtheorem{coro}[theo]{Corollary}
\newtheorem{lemma}[theo]{Lemma}

\theoremstyle{remark}

\newtheorem*{remark}{Remark}

\newenvironment{axiomlist}[1][MMMM]
    {\begin{list}{}{\settowidth{\leftmargin}{#1}} \setlength{\labelwidth}{\leftmargin} }
    {\end{list} }

\newcommand{\R}{\ensuremath{\mathds{R}}}
\newcommand{\C}{\ensuremath{\mathds{C}}}

\newcommand{\N}{\ensuremath{\mathds{N}}}

\DeclareMathOperator{\cov}{cov}
\DeclareMathOperator{\non}{non}

\newcommand{\Null}{\cal{N}}
\newcommand{\Meager}{\cal{M}}
\newcommand{\Ksigma}{\ensuremath{\cal{K}_{\sigma}}}

\DeclareMathOperator{\Pow}{Pow}
\DeclareMathOperator{\oc}{oc}

\newcommand{\cal}[1]{\ensuremath{\mathcal{#1}}}

\renewcommand{\frak}[1]{\ensuremath{\mathfrak{#1}}}
\newcommand{\thsym}[1]{\ensuremath{\mathsf{#1}}}

\newcommand{\set}[1]{\ensuremath{\left\{ #1 \right\}}}
\newcommand{\seqn}[1]{\ensuremath{\left( #1 \right)}}
\newcommand{\abs}[1]{\ensuremath{\left\lvert #1 \right\rvert}}

\newcommand{\expset}[2]{\ensuremath{{}^{#1}#2}}

\newcommand{\NonDecSeq}{\ensuremath{\expset{\N}{\N}{\uparrow}}}

\title{On the independence of a generalized statement of Egoroff's theorem from \thsym{ZFC}, after T.~Weiss}
\author[\sc R.~Pinciroli]{\sc Roberto Pinciroli}
\address{Scuola Normale Superiore\\56100 Pisa\\Italy}
\email{r.pinciroli@sns.it}
\thanks{I'm grateful to prof. T.~Weiss and prof. A.~Louveau for reading a preliminary version of this paper and suggesting some interesting questions for further study; I'd also like to thank dr. A.~Saracco and prof. G.~Tomassini, who informed me of the results in \cite{BSSW} and \cite{Wei}.}

\begin{document}

\begin{abstract}
We consider a generalized version (\thsym{GES}) of the well-known Severini\textendash Egoroff theorem in real analysis, first shown to be undecidable in \thsym{ZFC} by Tomasz Weiss in \cite{Wei}. This independence is easily derived from suitable hypotheses on some cardinal characteristics of the continuum like \frak{b} and \frak{o}, the latter being the least cardinality of a subset of $[0,1]$ having full outer measure.
\end{abstract}

\maketitle

In this paper we will consider the following \emph{Generalized Egoroff Statement}, which is a version ``without regularity assumptions'' of the well-known Severini\textendash Egoroff theorem from real analysis:

\begin{axiomlist}
\item[\thsym{GES}\ ] Given a sequence \seqn{f_n\ :\ n\in\N} of arbitrary functions $[0,1]\to\R$ converging pointwise to $0$, for each $\eta>0$ there is a subset $A\subseteq[0,1]$ of outer measure $\mu^*(A)>1-\eta$ such that \seqn{f_n} converges uniformly on $A$.
\end{axiomlist}

This conjecture first emerged from some questions about the behaviour of bounded harmonic functions on the unit disc in $\C$; in particular, it has been used in \cite{BSSW} to show the independence from \thsym{ZFC} of a strong Littlewood-type statement about tangential approaches.

Notice that in \thsym{GES} it is necessary to consider Lebesgue \emph{outer} measure to avoid simple counterexamples in \thsym{ZFC}:
\begin{prop}
There is a decreasing sequence \seqn{f_n\ :\ n\in\N} of functions $[0,1]\to\R$, converging pointwise to zero, such that every subset $A\subseteq[0,1]$ on which \seqn{f_n} converges uniformly has Lebesgue inner measure zero.
\end{prop}
\proof
By a theorem of Lusin and Sierpi\'nski there exists a partition of $[0,1]$ into countably many (in fact, even continuum many) pieces \set{B_n\ :\ n\in\N} each having full outer measure. Consider then the sequence \seqn{f_n} where, for every $n\in\N$, $f_n$ is the characteristic function of the subset $B_{\geq n}=\bigcup_{k\geq n} B_k$ of the unit interval: clearly \seqn{f_n(x)} converges monotonically to zero on every point $x\in[0,1]$; if \seqn{f_n} converges uniformly on a subset $A$, $A$ has to be disjoint from $B_{\geq\bar{n}}$ for some $\bar{n}\in\N$, so $\mu_{*}(A)\leq 1-\mu^*(B_{\geq\bar{n}})=0$.
\qed

Fix once and for all a decreasing vanishing sequence $\varepsilon=\seqn{\varepsilon_n}_{n\in\N}$ of positive real numbers, e.g. $\varepsilon_n=2^{-n}$; consider the following function, mapping a sequence of reals to its ($\varepsilon$-)\emph{order of convergence} to zero:
\begin{equation}\label{eq:OrderConvergence}
\begin{aligned}
&\oc\ :\ c_0\ \to\ \NonDecSeq,\quad\text{defined on each}\ a=\seqn{a_n}\in c_0\ \text{as} \\
&(\oc a)_n\;=\;\min\set{m\ :\ \forall l\geq m\ \big(\abs{a_l}\leq\varepsilon_n\big)},
\end{aligned}
\end{equation}
where $c_0$ denotes the set of infinitesimal real-valued sequences and $\NonDecSeq\subseteq\expset{\N}{\N}$ is the set of nondecreasing sequences of natural numbers. 

Using the natural identification of $\expset{\N}{(\expset{X}{\R})}$ with $\expset{X}{(\expset{\N}{\R})}$, we can view a sequence of real-valued functions $X\to\R$ converging pointwise to zero as a single function $F:X\to c_0$, and then study the associated order of convergence, $\oc F=\oc\circ F:X\to\NonDecSeq$:

\begin{lemma}\label{lemma:UniformConvergenceCompactOrders}
$F$ converges uniformly to zero if and only if the range of $\oc F$ is bounded in $(\expset{\N}{\N},\leq)$, where $\leq$ is the partial order of everywhere domination: $\alpha\leq\beta$ iff $\forall n\ (\alpha_n\leq\beta_n)$.
\end{lemma}
\proof
This is just a restatement of the definition of uniform convergence:
\begin{align*}
&F\ \text{converges uniformly to}\ 0\quad\leftrightarrow\\
&\leftrightarrow\quad\forall n\ \exists m\ \forall x\in X\ \forall l\geq m\ \big(\abs{F_l(x)}\leq\varepsilon_n\big)\quad\leftrightarrow\\
&\leftrightarrow\quad\exists \seqn{m_n}\in\expset{\N}{\N}\ \forall n\ \forall x\in X\ \big((\oc F(x))_n\leq m_n\big). \qedhere
\end{align*}

\begin{lemma}\label{lemma:ConstructionFunctionsGivenOrder}
For all $\varphi:X\to\NonDecSeq$ there exists a sequence $F$ of real-valued functions on $X$ converging pointwise to $0$ with order $\oc F=\varphi$.
\end{lemma}
\proof
It is sufficient to prove the lemma pointwise: given a nondecreasing sequence of natural numbers $\alpha\in\NonDecSeq$, we construct a sequence $a\in c_0$ converging to $0$ with order $\alpha$.
For that, just let
\begin{equation*}
a=\seqn{a_n}_{n\in\N} \quad\text{where}\quad a_n\;=\;\inf\set{\varepsilon_k\ :\ \alpha_k\leq n};
\end{equation*}
it is straightforward to check that this works, i.e. $\oc a=\alpha$.
\qed

Let $\mu^*$ be an upward continuous outer measure on a set $X$, i.e. an outer measure $\Pow X\to[0,+\infty]$ satisfying 
\begin{equation*}
A\ =\ \bigcup_{n\in\N} A_n \quad\rightarrow\quad \mu^*(A)=\lim_{n\to\infty}\mu^*\Big(\bigcup_{k<n}A_k\Big).
\end{equation*}
For every sequence $F$ of real-valued functions on $X$ converging pointwise to zero, consider the statement
\begin{axiomlist}[MMMMMMMM]
\item[$\thsym{GES}(X,\mu^*,F)$\ ] for each $M<\mu^*(X)$ there is a subset $A\subseteq X$ such that $\mu^*(A)>M$ and $F$ converges uniformly on $A$;
\end{axiomlist}
the Generalized Egoroff Statement relative to the space $(X,\mu^*)$ is the formula $$\thsym{GES}(X,\mu^*)=\forall F\ \thsym{GES}(X,\mu^*,F);$$ clearly our original statement $\thsym{GES}$ is just $\thsym{GES}([0,1],m^*)$, where $m^*$ is Le\-besgue outer measure on the unit interval $[0,1]\subseteq\R$.
Denote by \Ksigma\ the $\sigma$-ideal generated by the bounded subsets of $(\expset{\N}{\N},\leq)$; equivalently, \Ksigma\ is the family of those subsets which are bounded with respect to the order $\leq^*$ of eventual domination,
\begin{equation*}
\alpha\leq^* \beta \;\leftrightarrow\; \forall^{\infty}n\ (\alpha_n\leq\beta_n) \;\leftrightarrow\; \exists n\ \forall k\geq n\ (\alpha_k\leq \beta_k) \;\quad (\alpha,\beta\in\expset{\N}{\N}),
\end{equation*}
and \Ksigma\ is also the $\sigma$-ideal generated by the compact subsets of the Baire space \expset{\N}{\N} (see \cite{Bla}).

\begin{lemma}\label{lemma:CriterionGESSingleFamily}
$\thsym{GES}(X,\mu^*,F)$ holds iff there is a subset $Y\subseteq X$ of full outer measure (i.e. $\mu^*(Y)=\mu^*(X)$) such that $\oc F[Y]\in\Ksigma$.
\end{lemma}
\proof
Fix an increasing sequence of positive real numbers \seqn{M_n} with limit $\mu^*(X)$.
Assume $\thsym{GES}(X,\mu^*,F)$: by lemma \ref{lemma:UniformConvergenceCompactOrders}, for every $n\in\N$ there is a subset $A_n\subseteq X$ such that $\mu^*(A_n)>M_n$ and $\oc F[A_n]$ is bounded in \expset{\N}{\N}; taking $Y=\bigcup_{n\in\N}A_n$, $Y$ has full outer measure and $\oc F[Y]=\bigcup_{n\in\N}\oc F[A_n]$ is $\sigma$-bounded, as required.
Conversely, suppose that $\mu^*(Y)=\mu^*(X)$ and $\oc F[Y]\subseteq\bigcup_{n\in\N} B_n$, where each $B_n$ is a bounded subset of $(\expset{\N}{\N},\leq)$, and put
\begin{equation*}
A_n=(\oc F)^{-1}[B_0\cup\ldots\cup B_{n-1}]:
\end{equation*}
since $\oc F[A_n]$ is bounded, $F$ converges uniformly on every $A_n$ (lemma \ref{lemma:UniformConvergenceCompactOrders}); moreover, as $\mu^*$ is continuous and $Y\subseteq\bigcup_{n\in\N} A_n$, for all $m$ there is some $n$ such that $\mu^*(A_n)>M_m$, that is, $\thsym{GES}(X,\mu^*,F)$ holds.
\qed

\begin{theo}\label{theo:CriterionGES}
$\thsym{GES}(X,\mu^*)$ holds if and only if for all functions $\varphi:X\to\expset{\N}{\N}$ there is a subset $Y\subseteq X$ of full outer measure such that $\varphi[Y]\in\Ksigma$.
\end{theo}
This theorem provides a translation of \thsym{GES} into a purely set-theoretical statement.
\proof
The ``if'' direction follows directly from lemma \ref{lemma:CriterionGESSingleFamily} using $\varphi=\oc F$.
For the converse, consider the function $\varTheta$ which maps a sequence $\alpha=\seqn{\alpha_n}_{n\in\N}$ to the nondecreasing sequence $\seqn{\sum_{k\leq n} \alpha_k}_{n\in\N}$: it is a bijective order morphism $(\expset{\N}{\N},\leq)\to(\NonDecSeq,\leq)$ satisfying $\alpha\leq\varTheta(\alpha)$, therefore, for all $Y\subseteq\expset{\N}{\N}$, $\varTheta[Y]$ is ($\sigma$-)bounded iff $Y$ is ($\sigma$-)bounded. Assume $\thsym{GES}(X,\mu^*)$ and let $\varphi$ be a function $X\to\expset{\N}{\N}$; by lemma \ref{lemma:ConstructionFunctionsGivenOrder} there exists a sequence $F$ of real-valued functions converging pointwise to $0$ with $\oc F=\varTheta\circ\varphi$, so there is a set $Y\subseteq X$ of full outer measure such that $\varTheta[\varphi[Y]]=\oc F[Y]\in\Ksigma$ (lemma \ref{lemma:CriterionGESSingleFamily}), i.e. $\varphi[Y]\in\Ksigma$ as desired.
\qed

\begin{remark}
Theorem \ref{theo:CriterionGES} is still valid for measure spaces $(X,\mu)$ and the classical Egoroff Statement, provided that we only consider measurable maps $\varphi$ and measurable subsets $Y\subseteq X$. Thus theorem \ref{theo:CriterionGES} entails the Severini\textendash Egoroff theorem: if $\mu$ is finite and $\varphi:X\to\expset{\N}{\N}$ is measurable, the image measure $\varphi_{*}\mu$ is a finite Borel measure on \expset{\N}{\N}, hence it is regular and it is always supported by a $\sigma$-compact subset.
\end{remark}

Recall that the \emph{bounding number} $\frak{b}=\non(\Ksigma)$ (see \cite{Bla}) is the smallest possible size of a subset of \expset{\N}{\N} not belonging to \Ksigma.
We also denote with $\frak{o}(X,\mu^*)$ the least cardinality of a subset of $X$ having full outer measure and let $\frak{o}=\frak{o}([0,1],m^*)$%
\footnote{We haven't been able to find any specific name for this cardinal in the literature.}.

\begin{coro}\label{coro:ConsistencyGES}
Assuming $\frak{o}(X,\mu^*)<\frak{b}$, $\thsym{GES}(X,\mu^*)$ holds. In particular, $\frak{o}<\frak{b}$ implies \thsym{GES}%
\footnote{The latter fact has been pointed out by T.~Weiss and I.~Rec\l aw (see \cite{Wei}).}.
\end{coro}
\proof
Fix a subset $Y\subseteq X$ of full outer measure with $\abs{Y}=\frak{o}(X,\mu^*)$; then every function $\varphi:X\to\expset{\N}{\N}$ maps $Y$ onto a set of cardinality less than \frak{b}, hence $\varphi[Y]\in\Ksigma$.
\qed

We can also invoke theorem \ref{theo:CriterionGES} to prove sufficient conditions for the failure of \thsym{GES}. Precisely, we infer $\thsym{\neg GES}(X,\mu^*)$ by constructing (under suitable hypotheses) a set $Z\subseteq\expset{\N}{\N}$ of cardinality $\abs{Z}\geq\abs{X}$ such that all subsets of $Z$ belonging to \Ksigma\ have size less than $\frak{o}(X,\mu^*)$: once this is achieved, if $\varphi$ is any injection $X\to Z$, no subset $Y\subseteq X$ of full measure can be mapped onto an element of \Ksigma, because $\abs{\varphi[Y]}=\abs{Y}\geq\frak{o}(X,\mu^*)$. In order to state the next proposition, we recall that the \emph{dominating number} $\frak{d}\geq\frak{b}$ is the least cardinality of a cofinal subset of $(\expset{\N}{\N},\leq^*)$ and that a \emph{$\kappa$-Lusin set} is a subset $L\subseteq\R$ of cardinality $\kappa$ whose meager (i.e. Baire first category) subsets have size less than $\kappa$.

\begin{prop}\label{prop:ConsistencyNonGES}
Assume $\frak{o}(X,\mu^*)=\abs{X}=\kappa$; then $\thsym{GES}(X,\mu^*)$ fails in each of the following cases:
\begin{enumerate}
\item $\kappa=\frak{b}$;
\item $\kappa=\frak{d}$;
\item there exists a $\kappa$-Lusin set.
\end{enumerate}
\end{prop}
\proof
Following the plan outlined before stating the proposition, we try to build a ``$\kappa$-Lusin set'' $Z$ for the ideal \Ksigma\ instead of the ideal of meager sets. This is automatic under hypothesis (3): every (true) $\kappa$-Lusin set has the required properties, since all compact subsets of \expset{\N}{\N} have empty interior and thus every $\Ksigma$ set is meager.

Assume $\kappa=\frak{b}$ and let $\set{\alpha^{\xi}}_{\xi<\frak{b}}$ be an unbounded family in $(\expset{\N}{\N},\leq^*)$. By transfinite recursion we build a wellordered unbounded chain $Z=\set{\beta^{\xi}}_{\xi<\frak{b}}$ of length \frak{b}: after the construction of all $\beta^{\eta}$ for $\eta<\xi$, pick $\beta^{\xi}$ among the strict $\leq^*$-upper bounds of the set $\set{\alpha^{\xi}}\cup\set{\beta^{\eta}}_{\eta<\xi}$ (which has size less than \frak{b} and thus is $\leq^*$-bounded). It is clear that no $\leq^*$-bounded subset of $Z$ can be cofinal in $Z$, hence all \Ksigma\ subsets of $Z$ have cardinality $<\frak{b}$.

Finally, suppose $\kappa=\frak{d}$ and let $\set{\alpha^{\xi}}_{\xi<\frak{d}}$ be a cofinal family in $(\expset{\N}{\N},\leq^*)$. We build a set $Z=\set{\beta^{\xi}}_{\xi<\frak{d}}$ of cardinality \frak{d} by transfinite recursion as follows: after the construction of all $\beta^{\eta}$ for $\eta<\xi$, pick an element $\beta^{\xi}$ which is not $\leq^*$ any element of the set $\set{\alpha^{\eta}}_{\eta\leq\xi}\cup\set{\beta^{\eta}}_{\eta<\xi}$ (which has size less than \frak{d} and thus is not $\leq^*$-cofinal). $Z$ has the desired properties: $\seqn{\beta^{\xi}}_{\xi<\frak{d}}$ is a sequence without repetitions, hence $\abs{Z}=\frak{d}$, and moreover, if $A\subseteq Z$ is in $\Ksigma$, some $\alpha^{\xi}$ has to eventually dominate all elements of $A$, which implies that $A\subseteq\set{\beta^{\eta}}_{\eta<\xi}$ has cardinality less than \frak{d}.
\qed

\begin{coro}\label{coro:ConsistencyNonGES}
\thsym{GES} fails whenever at least one of the following hypotheses is satisfied:
\begin{enumerate}
\item $\frak{o}=\frak{d}=\frak{c}$ (the cardinality of the continuum);
\item there exists a \frak{c}-Lusin set and $\frak{o}=\frak{c}$;
\item there exists a \frak{c}-Lusin set and \frak{c} is a regular cardinal.
\end{enumerate}
\end{coro}
The last two conditions provide an affirmative answer (at least when \frak{c} is regular or it coincides with \frak{o}) to a question posed by T.~Weiss; he also noticed that there are models of \thsym{ZFC} (e.g. the iterated Mathias real model, where $\frak{o}=\frak{d}=\frak{c}$) which contain no \frak{c}-Lusin sets but nevertheless satisfy \thsym{\neg GES}.
\proof
Assumptions (1) and (2) are just particular instances of cases (2) and (3) respectively of proposition \ref{prop:ConsistencyNonGES}.
Moreover, hypothesis (3) is stronger than both (1) and (2): if $\kappa$ is a regular cardinal and there is a $\kappa$-Lusin set, then $\cov(\Meager)\geq\kappa$ and thus $\frak{d}\geq\cov(\Meager)\geq\kappa$ and $\frak{o}\geq\non(\Null)\geq\cov(\Meager)\geq\kappa$ (see \cite{BJ} for the relevant definitions of these cardinal characteristics associated to the $\sigma$-ideals \cal{M} of meager sets and \cal{N} of Lebesgue nullsets, as well as for the proofs in \thsym{ZFC} of the stated inequalities).
\qed

\begin{coro}[T.~Weiss]\label{coro:IndependenceGES}
\thsym{GES} is undecidable in \thsym{ZFC}.
\end{coro}
\proof
The hypothesis of corollary \ref{coro:ConsistencyGES}, and therefore \thsym{GES}, hold in the iterated Laver real model (see \cite{BJ} and the proof of theorem 1 in \cite{Wei}). On the other hand, $\frak{o}=\frak{d}=\frak{c}$ is certainly true (thus \thsym{\neg GES} holds) under the Continuum Hypothesis \thsym{CH} or just Martin's Axiom \thsym{MA}, which are consistent with \thsym{ZFC}.
\qed



\end{document}